\documentclass[12pt,letter]{article}
\usepackage{amsmath,amssymb,amsthm}
\usepackage{enumerate}

\usepackage{geometry}
\usepackage{euler}

\theoremstyle{plain}
\newtheorem{theo}{Theorem}
\newtheorem{lem}{Lemma}

\theoremstyle{definition}
\newtheorem{df}{Definition}

\theoremstyle{remark}
\newtheorem{rmk}{Remark}

\allowdisplaybreaks

\numberwithin{equation}{section}

\def\R{{\mathbb R}}
\def\Z{{\mathbb Z}}
\def\N{{\mathbb N}}

\def\E{{\mathbf E}}
\def\P{{\mathbf P}}
\def\1{{\mathbf 1}}
\newcommand{\calB}{\mathcal{B}}
\newcommand{\calN}{\mathcal{N}}
\newcommand{\calP}{\mathcal{P}}
\newcommand{\calV}{\mathcal{V}}

\def\eps{\varepsilon}
\let\phi=\varphi
\def\qed{\hfill\rule{.2cm}{.2cm}}

\allowdisplaybreaks
\title{Dynamics of the supermarket model}
\author{I. MacPhee$^{1}$ 
\and M.V.~Menshikov$^{1}$ 
\and M.~Vachkovskaia$^{2}$}

\begin{document}
\maketitle
{\footnotesize
\noindent
$^1$University of Durham,
Department of Mathematical Sciences,
South Road, Durham DH1 3LE, UK.\\
E-mails: i.m.macphee@durham.ac.uk, Mikhail.Menshikov@durham.ac.uk

\noindent
$^2$ Department of Statistics,
Institute of Mathematics, Statistics and Scientific Computation,
University of Campinas--UNICAMP,
P.O.\ Box 6065, CEP 13083--970, Campinas, SP, Brazil\\
E-mail: marinav@ime.unicamp.br

}

\begin{abstract}
  We consider the long term behaviour of a Markov chain $\xi(t)$ on
  $\Z^N$ based on the $N$ station supermarket model. Different routing
  policies for the supermarket model give different Markov chains. We
  show that for a general class of local routing policies, \emph{join
    the least weighted queue} (JLW), the $N$ one-dimensional
  components $\xi_i(t)$ can be partitioned into disjoint clusters
  $C_k$.  Within each cluster $C_k$ the \emph{speed} of each component
  $\xi_j$ converges to a constant $V_k$ and under certain conditions
  $\xi$ is recurrent in shape on each cluster. To establish these
  results we have assembled methods from two distinct areas of
  mathematics, semi-martingale techniques used for showing stability
  of Markov chains together with the theory of optimal flows in
  networks. As corollaries to our main result we obtain the stability
  classification of the supermarket model under any JLW policy and can
  explicitly compute the $C_k$ and $V_k$ for any instance of the model
  and specific JLW policy.
\end{abstract}



\section{Introduction} 
We consider the long term behaviour of a Markov chain $\xi(t)$ based
on the supermarket model of queueing theory.  In this model there are
$N$ stations, each of which processes jobs which queue there. Jobs
depart the system after their service is completed. The interesting
feature is that the stations support a neighbourhood structure of
non-empty sets of stations. Job streams arrive at these neighbourhoods
and upon arrival each job must be routed to a station within its
neighbourhood. The choice of queue can depend upon the current queue
lengths. Policies which route jobs based only upon information about
the queues in their own neighbourhoods are called \emph{local}. There
is considerable interest in the difference in performance between
systems with multi-station neighbourhoods and those with isolated
stations (so no routing) motivated by the work of Mitzenmacher and
others, see for example \cite{MRS} and \cite{VDK}, with some
sophisticated asymptotic work by Luczak and co-authors in \cite{LMcD}
and other papers.  The most commonly studied example of a local policy
is \emph{join the shortest queue} (JSQ). We consider a generalisation
of JSQ where each station $j$ has a weight factor $w_j > 0$ and each
job joins a \emph{least weighted} queue (JLW) at a station within its
neighbourhood (so JSQ is the case where all $w_j$ are equal).

A simple Markov model, $X(t)$ say, of such a system has independent
Poisson arrival streams to the neighbourhoods, exponential service
times at each station and lives on $\Z_+^N$. Our Markov chain $\xi(t)$
is based on $X(t)$ but we drop the requirement that the process is
non-negative so $\xi(t)$ lives on $\Z^N$ (we describe the transition
law in detail below). This enables us to exhibit behaviour of the
process that will only be seen for the queueing model $X(t)$ in large
deviation situations. While the Markov assumptions are strong we allow
general neighbourhood structures, arrival rates and service rates and
our results are about the long term behaviour of finite systems, not
large system asymptotics.

Our main result, Theorem \ref{speeds}, says that JLW policies induce
dependence between components $\xi_j$. Disjoint clusters of stations
appear (distinct from but determined by the neighbourhood structure
together with event rates) and at each station $j$ within cluster
$C_k$ say the drift rate $w_j\E \bigl(\xi_j(t+1) - \xi_j(t) \mid
\xi(t) \bigr) \to V_k$ for some constants $V_k$. It follows that the
weighted components $w_j\xi_j$ within a cluster are eventually much
closer to each other than to those in other clusters. Under some
constraints on the neighbourhoods and event rates we show in Theorem
\ref{shape} that the weighted process $w\xi(t)$ restricted to a
cluster is \emph{recurrent in shape}, an idea which appeared in Andjel
et al \cite{AMS} with some further application in \cite{MSV}. This
behaviour is caused by the routing policy and is akin to state space
collapse as discussed in several queueing network papers studying
heavy traffic e.g.\ Bramson \cite{B} and Kelly \& Williams \cite{KKLW}
though the time scales and techniques involved are entirely distinct.
Our results are established with semi-martingale/Lyapunov function
methods after the analysis of a carefully chosen deterministic flow
model on a graph.

This preliminary work on flows also leads to two new results for the
queueing model $X(t)$. Label the clusters and their drift rates so
that $V_1$ is the largest such rate. We show in Theorem
\ref{stability} that $X(t)$ is stable when $V_1 < 0$ and transient
when $V_1 > 0$. This result was shown for a system with identical
servers in a single neighbourhood by Weber \cite{W}, then for a
Markov system under JSQ by Foley and MacDonald \cite{FMcD} and then
for a system with more general arrival streams and service times again
under JSQ by Dai et al \cite{DHK} but we are not sure if it is known
for JLW. Writing $V_1(w)$ to indicate dependence upon the JLW weights
we have also shown in Theorem \ref{local}  that if $V_1(w) < 0$ for
some positive weights $w$ then $V_1(w') < 0$ for any set of positive
weights $w'$. In particular if $X(t)$ is stable under any version of
JLW it is also stable under JSQ. 

To establish these results we compare the behaviour of $\xi$ under JLW
to that under carefully chosen static policies which also cluster the
stations so that the weighted drift rates within clusters are
constant. In Theorem \ref{clusters} we show that the cluster structure
and drift rates $V_k$ can be determined for any neighbourhoods and
event rates by solving a particular flow problem on a bipartite graph.

\subsection{Model details and notation}
We will mostly consider two classes of simple Markov routing policies
described below. All jobs are of a single type but the servers have
different rates.  We assume here that service times at each station
$j$ are exponentially distributed with rate $\mu_j$ and that they are
independent of arrivals and other service times. Each job leaves the
system after completion of its service. We make no specific
assumptions about the queue discipline as we will not discuss waiting
times of individual jobs but we do assume the servers are non-idling
so when there are jobs in the queue at station $j$ the departure
process is of rate $\mu_j$.

The stations support a neighbourhood structure of non-empty sets of
stations $S_i \in \calP(C_0)$, the collection of subsets of $C_0 =
\{1,\,2, \ldots, N\}$. Jobs arrive at the neighbourhoods as
independent Poisson processes with rate $\lambda_i \geq 0$ at $S_i$
for each $i = 1$, $2$, $\ldots\,$. We allow some $\lambda_i = 0$ and
denote by $\calN(C_0)$ the neighbourhoods $S_i$ with $\lambda_i > 0$.
For simplicity we will usually write $i \in \calN(C_0)$ when we mean
$S_i \in \calN(C_0)$. To eliminate some trivial situations we suppose
the bipartite graph $G$, with nodes $\calN(C_0) \cup C_0$ and edges $E
= \{ (S_i, j) : j \in S_i\, \}$, is connected which ensures that model
cannot be trivially decomposed into independent components.

We are interested in the behaviour of the queue length process $X(t)$
on state space $\N^N$ ($\N$ denotes the non-negative integers)
and a related process, the random walk $\xi(t)$ with the same jump
rates as $X$ at positive states but not reflected at $0$ and hence
with state space $\Z^N$. The exact details of the jump rates depend
upon the routing policy so we discuss these now.

Upon arrival at neighbourhood $S_i$ a job is routed to a station
$j \in S_i$ where $j$ is chosen by some routing policy. We will mostly
consider two classes of simple Markov routing policies described
next. 

Define $\Delta_0 = \{ p \in [0, 1]^N : \sum_j p_j = 1 \}$ and for
each $i=1,\ \dots\,, |\calN(C_0)|$ let 
\[ \Delta_i = \{ p \in \Delta_0 : p_j = 0\ \text{for } j \notin
S_i \, \} \]
denote the unit simplex on coordinates $j \in S_i$. 

\begin{df}
A \emph{stationary Markov routing policy} is a mapping 
\[ \pi : \Z^N \times \{1, \, 2, \ldots\, , |\calN(C_0)| \} \to
  \Delta_0 \quad\text{such that } \pi(x, i) \in \Delta_i \,.  
\]
Under policy $\pi$ a job arriving at neighbourhood $S_i$ when the
process state is $x$ is routed to station $j \in S_i$ with probability
$\pi(x, i)_j$. We denote the space of routing policies by $\Pi$.

If $\pi$ does not depend upon $x$ we say it is \emph{static} and write
$\pi(i)$ for the routing distribution of arrivals at $S_i$. We denote
the space of static routing policies by $\Pi_{\text{stat}}$. \quad${}_\Box$
\end{df}

\begin{df}[Local routing policies]  \label{JLW}
Fix a set of positive weights $\{ w_j : j \in C_0, w_j > 0 \}$. 
For $x \in \Z^N$ let $\underline{wx}_i = \min_{l \in S_i} w_l x_l$
%
%
and let $\calB_i(x) = \{ j \in S_i : w_j x_j =  \underline{wx}_i\,\}$
denote the set of stations in $S_i$ with minimal weighted state value. The
\emph{join the least weighted} queue (JLW) routing policy is defined by
\[  \pi(x, i)_j = \left\{\begin{array}{ll} 
     1/|\calB_i(x)| &\text{ if } j \in \calB_i(x),  \\
     0          &\text{otherwise}.  \qquad {}_\Box
\end{array}  \right.
\]
\end{df}

\begin{rmk}
  The JLW policies are stationary Markov but not static. They are of
  practical interest as they are local and relatively simple to
  implement (only the queue lengths in an arrival's neighbourhood are
  needed to make its routing decision). There are several varieties of
  JLW (making different choices when $|\calB_i(x)| \geq 2$) and close
  variants like routing to stations where $w_j(x_j + 1)$ is minimal
  but the system behaviour at the level considered here is 
  much the same for all variants.
  
  Two particular cases have been studied for a variety of models.
  With $w_j = 1$ for each $j$ the policy is \emph{join the shortest
    queue} (JSQ). With weights $w_j = 1/\mu_j$ the policy
  is \emph{join the smallest workload} (JSW).  Another plausible
  choice of $w$ is to give most stations weight $1$ but protect some
  stations by making their $w_j$ larger.  \qquad${}_\Box$
\end{rmk}

The jumps and jump rates of the queue length process $X$ and related
random walk $\xi$ under stationary Markov routing policy $\pi$ are as
follows. The possible jumps change the current state $x$ by $\pm
e_j$, the unit vector in $\R^N$ with value $1$ in component $j$. Both
processes make up-jumps (a job arrives and is routed to a station)
with the same rates from all states $x \in \Z^N$ i.e.\ we have $x
\mapsto x + e_j$ at rate $\sum_{i : j \in S_i} \lambda_i \pi(x,
i)_j$. Downward jumps (job completions) $x \mapsto x - e_j$ occur at
rate $\mu_j$ at all states $x \in \Z^N$ for the random walk and
at all states $x \in \N^N$ with $x_j \geq 1$ for the queue length
process $X$. 

Under routing policy $\pi$ the drift rate $V$ of $\xi_j(t)$ at every
state $x$ and of $X_j(t)$ at $x$ with $x_j > 0$ satisfies
\begin{equation}  \label{V_j}
  V(j; x, \pi) = \sum_{i\in \calN(C_0)} \lambda_i \pi(x, i)_j -
\mu_j   
\end{equation}
which simplifies to $V(j; \pi) = \sum_{i\in \calN(C_0)} \lambda_i
\pi(i)_j - \mu_j$ for static policies. 

With our interest in JLW policies it is convenient to work with
weighted versions of $X$ and $\xi$. Let $w = (w_1\,, \ldots\, , w_N)$
be a set of positive station dependent weights. We will write $wx$ for
the vector $(w_1x_1, \ldots\,, w_Nx_N)$ for each $x \in \Z^N$ and with
this same convention write $wX$ and $w\xi$ for the weighted
processes. 

\begin{rmk} \label{weights}
  For any set of positive weights $w$ and under any policy
  $\pi$ the processes $X$ and $wX$ can be coupled so both reach state
  $0$ at the same times and hence both processes are recurrent or
  transient together under any fixed policies. We discuss stability of
  $X$ under different local policies below.  \qquad${}_\Box$
\end{rmk}

\begin{rmk} \emph{Neighbourhood-station interaction.} \quad
We have restricted our attention to station dependent service rates
$\mu_j$ here as models where jobs arriving at $S_i$ have service rates
$\mu_{ij}$ when routed to station $j$ show behaviour that is far from
optimal under local routing rules like JLW. 

Consider the following simple example. There are 3 stations $\{0, 1,
2\}$ and 3 neighbourhoods $S_i = \{i, [i+1] \}$, for $i = 0$, $1$, $2$
where $[i+1] = i+1 \pmod{3}$ and the service rates for each $i$ are
$\mu_{ij} = 1$ for $j = i$ but $\mu_{ij} = 1/2$ for $j = [i+1]$. The
arrival rates are $\lambda_i = 0.7$ for each $i$ and the policy that
routes all $S_i$ arrivals to station $i$ for each $i$ is clearly
stable and in fact minimizes the drift rate at each station. From the
symmetry of the situation JSQ (with ties split 50/50) 
sends half of all arrivals to the station where they receive the slow
service rate so the long run average service time is $(1 + 2)/2 = 1.5$
and as $0.7 \times 1.5 > 1$ the system will be unstable under
JSQ. \quad ${}_\Box$ 
\end{rmk}

\subsection{Results}
While our main interest is in local routing policies we start with
some results for static routing policies for the queue length process
$X$. To link them to the JLW policy we fix upon a set of positive
station dependent weights $w$. We are interested in static policies
that stabilize the system when this is possible i.e.\ policies that
keep the weighted drift rates (\ref{V_j}) small in the following
sense:
\begin{itemize}
 \item the maximal drift rate is as small as possible,
 \item the number of queues growing at maximal speed is minimal,
 \item the second largest drift rate is minimal on a minimal set of
   queues and so on.  
\end{itemize}

We will refer to any non-empty collection of stations $C
\subseteq C_0$ as a \emph{cluster} to separate it from association
with any particular arrival streams. For any cluster $C$
and any class of static policies $\Pi' \subset \Pi_{\text{stat}}$
define
\begin{equation} \label{generalV}
  \calV(C; w, \Pi') = \min_{\pi \in \Pi'} \max_{j\in C} w_j V(j ; \pi)
\end{equation}
i.e.\ the minimum (over policies in $\Pi'$) of the maximum drift rate
of weighted jobs over stations in $C$.  

\begin{theo}  \label{clusters} 
 Let $V_1 = \calV(C_0; w, \Pi_{\text{stat}})$.  We can decompose the set of
 stations into a hierarchy of disjoint clusters $C_1$, $\ldots\,$, $C_K$ 
 for some $1 \leq K \leq N$ with the following properties. 
\begin{enumerate}[(i)]
\item $C_1$ is the unique cluster $C$ such that $\calV(C; w,
  \Pi_{\text{stat}}) = V_1$ and $|C|$ is minimal.  
\item If $C_1 \neq C_0$ then for stages $k = 2$, $\ldots\,$ let 
\[ V_k = \calV(C_0 \setminus \cup_1^{k-1} C_n \,;\,w, \Pi_{k-1})   \]
 where $\Pi_{k-1}$ is the set of static policies that achieve $V_n$ on
 cluster $C_n$ for $n = 1$, $\ldots\,$, $k-1$. $C_k \subseteq C_0
 \setminus \cup_1^{k-1} C_n$ is the unique cluster that satisfies
 $\calV(C; w, \Pi_{k-1}) = V_k$ with minimal value of $|C|$. At each $k$ we
 have $V_k < V_{k-1}$. 

 For some $K \leq N$, $\cup_1^K C_i = C_0$ and the \emph{hierarchical
 minimax decomposition} is complete.  
\item $\Pi_K$ is non-empty. For any $\pi \in \Pi_K$ and each $j \in
  C_k$, $k = 1$, $\ldots\,$, $K$
\[  w_j \Bigl( \sum_{i\in \calN(C_0)} \lambda_i \pi(i)_j - \mu_j \Bigr)
  = V_k\,. \] 
\end{enumerate}
\end{theo}

It turns out that the clusters $C_i$ and drift values $V_i$ tell us a great
deal about the behaviour of the queue length process $X(t)$ and the
random walk model $\xi(t)$ under the JLW policy with the weights $w$
used in (\ref{generalV}). Under any static policy the queues at each
station are independent and for $\pi \in \Pi_K$ the drift rate at each
$j \in C_k$ is $V_k$. Under a JLW policy the queues will not be
independent in general. 

The next result concerns stability of $X(t)$ which for JSQ applied to
this Markov model first appears in Foley and McDonald \cite{FMcD} and
was extended by Dai et al \cite{DHK} for a model with non-Markov
arrival processes and non-exponential service.

\begin{theo}   \label{stability} 
  Choose a set of positive weights $w$. If $V_1 > 0$ then the queue
  length process $X(t)$ is transient under any policy. If $V_1 < 0$
  then the queue length process $X(t)$ is positive recurrent under the
  JLW policy with weights $w$ and under any $\pi \in \Pi_K$.
\end{theo}

While this stability condition is known in the JSQ case (all $w_j =
1$) our approach via the hierarchical minimax decomposition with
static policies allows us to compare the stability of $X$ under
different local routing policies.

\begin{theo}  \label{local}
  Let $V_1(w)$ denote the maximal drift rate obtained from the
  hierarchical minimax decomposition with weights $w$. If $V_1(w) < 0$
  for some positive $w$ then $V_1(w') < 0$ for all positive weights
  $w'$. Thus if $X$ is stable under any JLW policy it is stable under
  all such policies.
\end{theo}

Next we consider the dynamics of the random walk process $\xi$ under
local policies with any fixed choice of weights $w$. The JLW policy
makes the queue workloads within clusters dependent and our main
result is that $V_k$ can be interpreted as the rate of change to the
weighted queue length at each queue in cluster $C_k$ under JLW
routing. This is a new observation for this model. In stable cases
($V_1 < 0$) this behaviour of the queueing model $X(t)$ will only be
seen in large deviation situations due to the reflection of the
process at $0$. This is why we have introduced the random walk $\xi$.

\begin{theo}  \label{speeds}
  The random walk $\xi(t)$ under the JLW policy eventually displays
  the hierarchical minimax structure. Specifically the weighted random
  walk $w\xi(t)$ eventually has drift rate $V_k$ on each cluster
  $C_k$, that is for small enough $\eps > 0$ and any finite initial
  configuration $\xi(0) = x_0$ there exists a random time $t(\eps)$
  such that for each cluster $C_k$ and for any $t > t(\eps)$  
\[ \left| \frac{w_j\xi_j(t)}{t} - V_k \right| < t^{-\eps} \quad\text{for
  each } j \in C_k \ .\] 
\end{theo}

\begin{rmk}
  In cases where $V_k > 0$ for some clusters $C_k$ this result can be
  extended to apply to the queueing model $X(t)$ on such clusters.
  \qquad${}_\Box$ 
\end{rmk}

Under some slightly stronger conditions on the internal structure of
the clusters we find that the JLW policy causes $w \xi$ to exhibit
some remarkable behaviour which we call \emph{recurrence in shape}.

\begin{df}
  For any Markov process $Y(t)$ on $\Z^N$ let $Y_C(t) = (Y_j(t) : j
  \in C)$ denote the process (perhaps not Markov) on components in $C$
  for any $C \subset \{1, 2, \ldots\,, N\}$. 
  For any $j_0 \in C$ we say that $Y_C(t)$ is \emph{recurrent in
  shape} (or $Y(t)$ is recurrent in shape on $C$) when the process 
  $(Y_j(t) - Y_{j_0}(t)) : j \in C)$ is recurrent.  \quad${}_\Box$
\end{df}

This notion was applied to a single cluster storage model in
\cite{MSV}. To apply it here we need a slightly more refined
description of the internal structure of the clusters. For $\pi \in
\Pi_{\text{stat}}$ let $G(\pi)$ denote the bipartite graph with nodes
$\calN(C_0) \cup C_0$ and edges $E = \{ (S_i, j) : \pi(i)_j > 0 \}$. 

\begin{df}
  Cluster $C \subseteq C_K$ is \emph{bonded} if there exists a policy
  $\pi \in \Pi_k$ such that for every pair $j, m \in C$ there is a
  path in $G(\pi)$ from $j$ to $m$.  \qquad${}_\Box$
\end{df}

\begin{theo} \label{shape} 
  Under the JLW policy with weights $w$ the weighted random walk
  $w\xi(t)$ is \emph{recurrent in shape} on each bonded sub-cluster of
  $C_k$, $k = 1$, $\ldots\,$, $K$.
\end{theo}

The restriction to bonded clusters is necessary as in general a
cluster can split into two or more independent parts with the same
drift rates.

The clustering behaviour described here is probably also exhibited by
the supermarket model with more general arrival streams and
non-exponential service times and we are looking for Lyapunov
functions that will allow us to extend our arguments at the necessary
points. It is possible that similar behaviour will persist in similar
systems which have Jackson-style feedback though there are many
complications here, see for example Dai et al \cite{DHK}. In fact it
was in trying to understand \cite{DHK} that we discovered the results
described here.

\section{Proofs}
\subsection{Preliminaries for Theorem \ref{clusters}}
We now describe a flow based decomposition of the system based on some
ideas from the max-flow, min-cut theorem of Ford and Fulkerson for
models of flow in networks. Then we show that the hierarchical minimax
decomposition coincides with the flow based decomposition. In this
section we only consider static routing policies and we work with a
fixed set of positive weights $w$ as before. 

\paragraph{Flow based decomposition}  
Recall that any non-empty set of stations is called a
\emph{cluster}. For each cluster $C \subseteq C_0$ let $\calP(C)$ be
the collection of subsets of $C$ and $\calN(C) = \{i \in \calN(C_0) :
S_i \subset C\}$ the collection of neighbourhoods with $\lambda_i > 0$
supported by $C$. The system \emph{restricted} to $C$ consists of the
stations in $C$ with the arrival streams to neighbourhoods in
$\calN(C)$. Such a restriction can be achieved by applying a routing
policy from 
\[ \Pi(C) = \{ \pi \in \Pi_{\text{stat}} : \pi(i)_j = 0\ \text{if } i
    \notin \calN(C),\  j \in C \},
\]
the set of static policies which are \emph{consistent} with this
decomposition. 

We also introduce the idea of the system \emph{reduced onto} a cluster
$D$ by removing the stations in $C_0 \setminus D$ together with the
arrival streams to neighbourhoods in $\calN(C_0 \setminus D)$. For
each non-empty $S \in \calP(D)$ define $\sigma_S(D) = \{S_i \in
\calN(C_0) : S_i \cap D = S \}$ i.e.\ the collection of neighbourhoods
(if any) that coincide with $S$ on $D$, and merge all arrival
streams in $\sigma_S(D)$ to get one with rate 
\begin{equation}   \label{lambdaD}
  \lambda_S(D) = \sum_{i\in \sigma_S(D)} \lambda_i\ .  
\end{equation}
Now for $C \subseteq D$ let $\calN_D(C) = \{S \in \calP(C) :
\lambda_S(D) > 0 \}$.  Combining these two notions we see there are
static routing policies which act to decompose the original system
into one restricted to a cluster $C_1$ and an independent system
reduced onto $D = C_0 \setminus C_1$ which can be further decomposed
as desired. These steps can be repeated to sequentially decompose the system.

We must also consider the drift rates at stations under any such
decomposition. For the system \emph{reduced onto} $D$ define, for any
cluster $C \subseteq D$, the average \emph{$D$-reduced restricted
  drift} on $C$ by
\begin{equation}   \label{W_D}
  W_D(C) = \frac{1}{|C|} \left( \min_{\pi \in \Pi(C)} \sum_{j\in
  C} w_j  \Bigl( \sum_{i\in \calN_D(C)} \lambda_i(D)\pi(i)_j -
  \mu_j  \Bigr) \right) \ .  
\end{equation}
We will not change the weights $w = (w_j)$ during the decomposition so
we do not indicate $W_D$'s dependence upon $w$. 

\begin{rmk}
The conditions $W_D(C) \leq 0$ for every cluster $C \subset D$ are
closely related to the sufficient conditions for a matching on a
bipartite graph. In the special case where each $w_j = 1$ we can
reverse the order of summation above and sum out the dependence on
$\pi$ to get 
\[  W_D(C) = \frac{1}{|C|} \Bigl( \sum_{i\in \calN_D(C)} \lambda_i(D)
   - \sum_{j\in C} \mu_j  \Bigr) \ .
\]
Hence $W_D(C) \leq 0$ only when the total service rate of
servers in $C$ is at least as large as the total rate of arrivals to
neighbourhoods entirely supported by $C$. \qquad ${}_\Box$
\end{rmk}

Jobs arriving at neighbourhoods in $\calN(C)$ cannot be routed to
stations outside $C$ but it may be possible to route jobs from other
neighbourhoods to stations in $C$ so, by comparison with (\ref{V_j}),
for any policy $\pi$ at any $x$ where $\pi$ is consistent with
reduction onto $D$ and for any $C \subset D$ we have 
\begin{equation}  \label{Vbound}
  \sum_{j \in C} w_j V(j ; x, \pi) = \sum_{j \in C} w_j \Bigl( \sum_{i
  \in \calN_D(D)} \lambda_i(D) \pi(x, i)_j - \mu_j \Bigr) \geq  |C| W_D(C)
\end{equation}
i.e.\ the $D$-reduced restricted drift is a lower bound for the
average drift over $C$ under any policy $\pi$ on the $D$-reduced
system and this lower bound is reached by some policies.   

There is one further small result which it is convenient to separate
out from the proof of Theorem \ref{clusters}. 
\begin{lem} \label{W(AB)}
  Suppose $A$, $B$ are clusters contained within $D$ such that 
\[ W_D(A) = W_D(B) = \max_{C \subset D} W_D(C) = v\,.
\]
  Let $\calN_1 = \calN_D(A \cup B) \setminus (\calN_D(A) \cup
  \calN_D(B))$. Then $W_D(A \cup B) = v$ and $\sum_{\calN_1}
  \lambda_i(D) = 0$.  
\end{lem} 

\paragraph{Proof} Let $H = A \cap B$ and $\calN_2 = \calN_D(B) \setminus
\calN_D(H)$ and note that $\calN_D(B \setminus H) \subset
\calN_2$. Let $\Pi(H)_v$ denote the set of policies that achieve
$\sum_{j\in H} w_j \bigl( \sum_{\calN_D(H)} \lambda_i(D)\pi(i)_j -
\mu_j \bigr) \leq |H|v$. By maximality of $v$ we have $W_D(H) \leq v$
and hence $\Pi(H)_v$ is non-empty and in fact it contains policies
that achieve the $D$-reduced restricted drifts for the clusters $A$ and
$B$.   

For any policy $\pi$ we have  
\[ |B|W_D(B) \leq \sum_{j\in H} w_j \sum_{\calN_D(H)} \lambda_i(D)
  \pi(i)_j + \sum_{j\in B} w_j \sum_{i \in \calN_2} \lambda_i(D)
  \pi(i)_j - \sum_{j \in B} w_j \mu_j \,. 
\]
For any policy $\pi \in \Pi(H)_v$\ it follows that 
\[ |B \setminus H|v \leq \sum_{j \in B} w_j \sum_{i \in \calN_2}
     \lambda_i(D)\pi(i)_j - \sum_{B \setminus H} w_j \mu_j  
\]
and from this we have, for any $\pi \in \Pi(H)_v$\,, 
\begin{eqnarray*}  
  \lefteqn{ \sum_{j \in A \cup B} w_j \Bigl( \sum_{i \in \calN_D(A \cup B)}
    \lambda_i(D) \pi(i)_j - \mu_j \Bigr)} \hspace*{3cm} \\
  & = &\sum_{j \in A} w_j \Bigl( \sum_{i \in\calN_D(A)} \lambda_i(D)\pi(i)_j
    - \mu_j \Bigr) + \sum_{j \in B} w_j \sum_{i \in
    \calN_2} \lambda_i(D)\pi(i)_j   \\
  & &\hspace*{3cm} {} - \sum_{j \in B \setminus H} w_j \mu_j + \sum_{j
    \in A \cup B} w_j \sum_{i \in \calN_1} \lambda_i(D)\pi(i)_j   \\ 
  & \geq &v(|A| + |B \setminus H|) + \sum_{j \in A \cup B} w_j
    \sum_{i \in \calN_1} \lambda_i(D)\pi(i)_j \,.   
\end{eqnarray*}
By maximality of $v$ and $|A \cup B| = |A| + |B \setminus
H|$ we see that $\sum_{\calN_1} \lambda_i(D) = 0$ and $W_D(A \cup B) =
v$ as required.  \quad\qed

\subsection{Proof of Theorem \ref{clusters}}
The routing schemes we consider here minimize the maximum drift by
directing some arrivals from heavily loaded parts of the network to
less loaded parts. Part (i) of this Theorem deals with the most
heavily loaded part first.

\noindent $(i)$ Using the flow based scheme we define at stage 1, $D =
C_0$ and 
\[  W_1 = \max_{C\subseteq D} W_D(C)\ , \quad 
   \overline{C}_1 = \cup \{C \subseteq D : W_D(C) =  W_1 \}
\]
and we now show that $W_D(\overline{C}_1) = W_1$. If there is only one
$C$ with $W_D(C) = W_1$  we are done so suppose there are more and
take any two distinct clusters $C$, $C'$ such that $W_D(C) = W_D(C') =
W_1$. By Lemma \ref{W(AB)}, $W_D(C \cup C') = W_1$ from which it soon
follows that $W_D(\overline{C}_1) = W_1$ and hence $\overline{C}_1$ is
the unique maximal cluster with restricted drift value $W_1$.

Next we must relate $W_1$ (an average drift on a cluster) to $V_1$ (a
maximal drift within a cluster). We start by showing that there
exist policies $\pi' \in \Pi(\overline{C}_1)$ such that $w_jV(j ;
\pi') = W_1$ for each $j \in \overline{C}_1$. If not then let
$\Pi(\overline{C}_1)^* \subset \Pi(\overline{C}_1)$ denote the
set of policies $\pi$ that achieve the restricted drift $W_1$ i.e.\ 
\[ \sum_{\overline{C}_1} w_j V(j ; \pi) = |\overline{C}_1| W_1 \ .
\]
Next pick $\hat\pi \in \Pi(\overline{C}_1)^*$ that
achieves $\hat{v} = \min_{\Pi(\overline{C}_1)^*} \max_{j \in
\overline{C}_1} w_jV(j ; \pi)$ at some stations in
$\overline{C}_1$. By assumption $\hat{v} > W_1$.  
Let $\hat{C} = \{j \in \overline{C}_1 : w_jV(j ; \hat\pi) =
\hat{v} \}$ and consider the restricted drift $W_D(\hat{C})$ on
$\hat{C}$. As $\hat{v}$ is minimal any calls that can be routed out of
$\hat{C}$ by $\hat{\pi}$ will be and so $\hat{v} = W_D(\hat{C})$ but
this implies $W_D(\hat{C}) > W_1$ so by maximality of
$W_1$ we must have $\hat{v} = W_1$. 

Hence there exist policies $\pi' \in \Pi(\overline{C}_1)^*$ such that
$w_jV(j ; \pi') \leq W_1$ for $j \in \overline{C}_1$. As $\sum_{j \in
  \overline{C}_1} w_jV(j ; \pi') = |\overline{C}_1| W_1$ for such
policies we must have $w_jV(j ; \pi') = W_1$ for each $j \in
\overline{C}_1$ as required. By maximality of $\overline{C}_1$ we
cannot have $w_jV(j ; \pi') = W_1$ for all $\pi'$ for any $j \notin
\overline{C}_1$. 

Now pick $\hat{\pi} \in \Pi(\overline{C}_1)$ such that $w_j V(j ;
\hat{\pi}) = W_1$ for all $j \in \overline{C}_1$ and $w_j
V(j ; \hat{\pi}) < W_1$ for $j \notin \overline{C}_1$. Then 
\[ V_1 = \min_{\pi \in \Pi_{\text{stat}}} \max_{j \in D} w_j V(j ;
  \pi) \leq \max_{j \in D} w_j V(j ; \hat{\pi}) = W_1 
\]
but from (\ref{Vbound}) we have, for any $\pi \in \Pi_{\text{stat}}$, 
\[ \max_{j \in D} w_j V(j ; \pi) \geq \frac{1}{|\overline{C}_1|}
\sum_{j \in \overline{C}_1} w_j V(j ; \pi) \geq W_1\ .  
\]
Hence $V_1 = W_1$ and no policy with maximal drift rate $V_1$ can
achieve drift rate less than $V_1$ at any station $j \in
\overline{C}_1$. As $w_j V(j ; \hat{\pi}) = V_1$ 
for $j \in \overline{C}_1$, $w_j V(j ; \hat{\pi}) < V_1$ for $j \notin 
\overline{C}_1$ it follows that $C_1 = \overline{C}_1$. 
This completes part (i) of the theorem. 

\noindent $(ii)\ \&\ (iii)$ For systems where $C_1 \neq C_0$ we can
continue for stages $k \geq 2$. At stage $k$ reduce the system onto
$D_k = C_0 \setminus \cup_1^{k-1} C_i$ and define
\[   W_k = \max_{C\subseteq D_k} W_{D_k}(C)\ , \quad 
    \overline{C}_k = \cup \{C \subseteq D_k : W_{D_k}(C) = W_k \}.
\]
At each stage $V_k = W_k$ and $C_k = \overline{C}_k$ follow as in
stage 1. That $V_k < V_{k-1}$ follows from the restricted drifts for
if $V_k \geq V_{k-1}$ then $W_{D_{k-1}}(C_{k-1} \cup C_k) \geq
W_{D_{k-1}}(C_{k-1})$ and $|C_{k-1} \cup C_k| > |C_{k-1}|$ in
contradiction to $|C_{k-1}|$ being maximal.

After a finite number of stages $K$ ($K \leq N$, the number of
stations) we will have $C_K = D_K$ which completes the
decomposition.  \qed

\begin{rmk}  \label{V_k} \emph{(conservation of mass)}
  Consider static policies $\pi$ that are consistent with reduction of the
  system onto $D = \cup_{l=k}^K C_l$ and route all flow possible out
  of $C_k$.  If under such a $\pi$ we have, for each $j \in C_k$, $w_j
  \bigl(\sum_{\calN_D(C_k)} \lambda_i(D) \pi(i)_j - \mu_j \bigr) = V$
  where $V$ is constant then $V = V_k$. To see this sum the equations
  $w_j\bigl(\sum_{i \in \calN_D(C_k)} \lambda_i(D) 
  \pi(i)_j - \mu_j \bigr) = V$ over $j \in C_k$ to get 
\begin{eqnarray*}
  V \sum_{j\in C_k} 1/w_j & = &\sum_{j \in C_k} \Bigl(\sum_{i \in
  \calN_D(C_k)} \lambda_i(D) \pi(i)_j - \mu_j \Bigr) \\
  & = &\sum_{\calN_D(C_k)} \sum_{C_k} \lambda_i(D) \pi(i)_j -
  \sum_{C_k} \mu_j = \sum_{\calN_D(C_k)} \lambda_i(D) - \sum_{C_k}
  \mu_j\ . 
\end{eqnarray*}
This equation does not depend upon $\pi$ and is satisfied by
$V_k$. We establish the same result for more general policies in the
proof of Theorem \ref{speeds} but there is no equivalent result for
models with service rates $\mu_{ij}$ that depend upon the routing
decision $S_i$ to $j$. \hspace*{1cm}${}_\Box$
\end{rmk}

\subsection{Proof of Theorem~\ref{stability}}
For any static routing policy $\pi$, the arrivals to stations $1,
\ldots, N$ are independent Poisson processes, arrivals to station $j$
having rate $\sum_{i \in \calN(C_0)} \lambda_i \pi(i)_j$. The system
is simply a collection of $N$ independent M/M/1 queues so it is
ergodic under $\pi \in \Pi_K$ if and only if all $V_k < 0$ and this is
implied by $V_1 < 0$.

Under local policies the queues become dependent and we use Lyapunov
or test function results to establish transience or recurrence
properties. We briefly state a couple of well known results that we use
a few times in what follows. 
\begin{theo}  \label{Foster}
Suppose $\{ X_n \}$ is an irreducible Markov chain on a countable state
space $\mathcal{S}$ and $f : \mathcal{S} \to \R^+$. Let $\Delta f_n =
f(X_{n+1}) - f(X_n)$. 
\begin{enumerate}[(i)]
\item If there are constants $c > 0$, $d > 0$ and $\eps > 0$ such that
  $|\Delta f_n| < d$ a.s.\ and \\
  $\E( \Delta f_n \mid X_n = x) > \eps$ for all $x \in \{x : f(x) > c
  \}$ then $\{ X_n \}$ is transient.  
\item If there is a constant $\eps > 0$ and a finite set $A \subset
  \mathcal{S}$ such that $\E( f(X_{n+1}) \mid X_n = x) < \infty$ for
  $x \in A$ and $\E(\Delta f_n \mid X_n = x) \leq -\eps$ for $x \in
  \mathcal{S} \setminus A$ then $\{ X_n \}$ is positive recurrent. 
\end{enumerate}
\end{theo}

\noindent\emph{Proof:} part (i) is a special case of Theorem 2.2.7 in
\cite{FMM} -- note the need for bounded jumps. Part (ii) is Foster's
criterion which can be found in many places, for instance Theorem
2.2.4 of \cite{FMM} or Proposition I.5.3 of \cite{Asm}.  \qed
\smallskip  

Now we return to the proof of Theorem \ref{stability}. If $V_1 > 0$ we
show transience under any policy $\pi$ by using the Lyapunov function
$L(x) = \sum_{j \in C_1} w_j x_j$ with $X(t)$'s jump chain $X_n$. The
rate of events for $X(t)$ at state $x$ is given by \( \alpha(x) =
\sum_{\calN(C_0)} \lambda_i + \sum_{C_0} \mu_j1_{\{x_j > 0\}} \) which
is bounded. Using the notation $\mathbf{0}$ for the state with every
$x_j = 0$, $\mathbf{1}$ for the state where every $x_j = 1$ we have
\begin{equation} \label{event_rate}
  \alpha(\mathbf{0}) = \sum_{\calN(C_0)} \lambda_i \leq \alpha(x)
  \leq \sum_{\calN(C_0)} \lambda_i + \sum_{C_0} \mu_j =
  \alpha(\mathbf{1}) \ .    
\end{equation}
Let $\Delta L_n = L(X_{n+1}) - L(X_n)$. We have $|\Delta L_n| \leq
\max_j w_j$ at all states of the system. For any policy $\pi$, inequality
(\ref{Vbound}) applied to cluster $C_1$ with $D = C_0$ implies
\[ \alpha(x) \E_{\pi}(\Delta L_n \mid X_n = x) = \sum_{j \in C_1} w_j
V(j ; x, \pi) \geq |C_1|V_1 > 0 
\]
at every state $x \in \N^N$. By Theorem \ref{Foster}(i) the jump chain
$X_n$ is transient under $\pi$ and hence so is the queue length
process $X(t)$.  

It remains to show that $X(t)$ is stable under the JLW policy when
$V_1 < 0$. As the event rates lie in the interval
$[\alpha(\mathbf{0}),  \alpha(\mathbf{1})]$ we can work with the jump
chain $X_n$ instead. It is convenient to work with a quadratic
Lyapunov function here. 

Let $q(x) = \frac{1}{2} x^TQx$ where $Q$ is a real, symmetric $N
\times N$ matrix and let $e_j$ denote the unit vector with $e_{jj} =
1$. For $\delta \in \{-1, 1\}$ we have 
\[ q(x + \delta e_j) - q(x) = \delta e_j^T Qx +  \frac{1}{2} Q_{jj} \] 
and we need to compute $\E_{\pi}(\Delta q_n \mid X_n =x)$ where
$\Delta q_n = q(X_{n+1}) - q(X_n) $. With indicator functions $D_j =
1_{\{\text{departure from station } j\}}$, $A_i = 1_{\{\text{arrival
at } S_i\}}$, $R_{ij} = 1_{\{S_i\ \text{arrival routed to } j\}}$  we
have 
\begin{eqnarray}
   \Delta q_n & = &\sum_{j \in C_0} \Bigl\{ D_j \Bigl[ q(X_n - e_j) -
    q(X_n) \Bigr] + \sum_{i\in \calN(C_0)} A_i R_{ij} \Bigl[ q(X_n + e_j) -
    q(X_n) \Bigr] \Bigr\} \nonumber\\
    & = &\sum_{j \in C_0} \Bigl(\sum_{i\in \calN(C_0)} A_iR_{ij} - D_j\Bigr)
    e_j^TQX_n + \frac{1}{2} \sum_{j \in C_0} Q_{jj} \Bigl( \sum_{i \in
    \calN(C_0)} A_i R_{ij}  + D_j \Bigr)\ . \label{delta_q}
\end{eqnarray}
 Also, for any policy $\pi$, we have 
\[  \alpha(x) \E_{\pi}(D_j \mid X_n = x) = \mu_j 1_{\{x_j > 0\}},  \qquad 
  \alpha(x) \E_{\pi}(A_i R_{ij} \mid X_n = x) = \lambda_i \pi(x, i)_j\ . 
\]
For the specific function $S(x) = \frac{1}{2}\sum_{j \in C_0} w_j
x_j^2$ we have $Q_{jj} = w_j$ and $e_j^TQx = w_j x_j$.  
Further, for $\pi \in \Pi_K$, $j \in C_k$ we have
\( w_j \Bigl( \sum_{\calN(C_0)} \lambda_i \pi(i)_j - 
\mu_j\Bigr) = V_k\) and so for $\pi \in \Pi_K$, 
\begin{equation}  \label{superMG}
  \alpha(x) \E_{\pi}(\Delta S_n \mid X_n = x) = \sum_k V_k \sum_{j
  \in C_k} x_j + \beta_{\pi}(x) 
\end{equation}
(any coefficient is OK for terms $x_j = 0$) where $\Delta S_n =
 S(X_{n+1}) - S(X_n)$ and
\[ \beta_{\pi}(x) =  \frac{1}{2} \sum_{j \in C_0} w_j\Bigl( \mu_j
 1_{\{x_j > 0\}} + \sum_{i \in \calN(C_0)} \lambda_i \pi(i)_j \Bigr)
 \leq \frac{1}{2} \alpha(\mathbf{1})   \max_j w_j \ . \]  
As each $V_k < 0$ the process $S(X_n)$ is a good supermartingale at
all but a finite subset of $\N^N$ under any $\pi \in \Pi_K$. 

We complete the proof by comparing the behaviour of $S(X_n)$ under JLW
with its behaviour under $\pi \in \Pi_K$. Let $\E_L$ denote
expectation under the JLW policy. Only the variables $R_{ij}$ are
controlled by the routing policy so by comparison with (\ref{delta_q}) 
\begin{eqnarray*}
 \lefteqn{\alpha(x) \Bigl( \E_L(\Delta S_n \mid X_n = x)  -
   \E_{\pi}(\Delta S_n \mid X_n = x) \Bigr)   } \qquad \\ 
 & = &\alpha(x) \sum_{j \in C_0} w_j (x_j + 1/2) \sum_{i \in \calN(C_0)}
   \Bigl[ \E_L(A_i R_{ij} \mid X_n = x) - \E_{\pi} (A_i R_{ij} \mid
   X_n = x) \Bigr]  \\
 & = &\alpha(x) \sum_{i \in \calN(C_0)} \sum_{j \in S_i} w_j x_j \Bigl[
 \E_L(A_i R_{ij} \mid X_n = x) - \E_{\pi} (A_i R_{ij} \mid X_n = x)
 \Bigr]  \\
 & &\qquad {} + \frac{\alpha(x) }{2}\sum_{j \in C_0} w_j \sum_{i \in
   \calN(C_0)} \Bigl[ \E_L(A_i R_{ij} \mid X_n = x) - \E_{\pi} (A_i
   R_{ij} \mid  X_n = x) \Bigr] 
\end{eqnarray*}
For each neighbourhood recall that $\underline{wx}_i = \min_{j \in
S_i} w_j x_j$ and let $\hat{w}_i = \max_{j \in S_i} w_j$. For the
second part of this sum we have the simple bound  
\begin{eqnarray*}
   \lefteqn{\alpha(x) \sum_{j \in C_0} w_j \sum_{i \in \calN(C_0)}
   \Bigl[ \E_L(A_i R_{ij} \mid X_n = x) - \E_{\pi} (A_i R_{ij} \mid
   X_n = x) \Bigr] } \hspace{3cm} \\
  & = &\sum_{\calN(C_0)} \sum_{j \in S_i} w_j
  \frac{\lambda_i}{|\calB_i(x)|} \, 1_{\{w_jx_j =
    \underline{wx}_i \}} - \sum_{C_0} w_j \sum_{\calN(C_0)} \lambda_i
  \pi(i)_j \\
  & \leq &\sum_{\calN(C_0)} \lambda_i \Bigl( \hat{w}_i - \sum_{j \in
   S_i} w_j \pi(i)_j \Bigr)
\end{eqnarray*}
 which does not depend on $x$. We now show the first part
 of the sum is negative. The JLW policy routes arrivals to stations in
 $\calB_i(x) \subset S_i$ where we have $w_jx_j =
 \underline{wx}_i$. Hence
\begin{eqnarray*}
\lefteqn{\alpha(x)  \sum_{\calN(C_0)} \sum_{j \in S_i} w_j x_j \Bigl[
 \E_L(A_i R_{ij} \mid X_n = x) - \E_{\pi} (A_i R_{ij} \mid X_n = x)
 \Bigr]} \qquad \\
 & = &\alpha(x)  \sum_{\calN(C_0)} \Bigl(\underline{wx}_i + (w_j x_j -
 \underline{wx}_i) \Bigr) \sum_{j \in S_i} \Bigl[  \E_L(A_i R_{ij}
 \mid X_n = x) - \E_{\pi} (A_i R_{ij} \mid X_n = x)  \Bigr] \\
 & = &\sum_{\calN(C_0)} \underline{wx}_i (\lambda_i - \lambda_i) 
 + \sum_{\calN(C_0)} \sum_{S_i \setminus \calB_i(x)} 
  (w_j x_j - \underline{wx}_i) (0 - \lambda_i \pi(i)_j) \leq 0
\end{eqnarray*}
with strict inequality except when the $w_jx_j$ are all equal or
$\pi(i)_j = 0$ outside $\calB_i(x)$ for every $S_i$. Combining this
with (\ref{superMG}) where $\pi \in \Pi_K$ we have 
\begin{eqnarray*}
  \alpha(x) \E_L (\Delta S_n \mid X_n = x) & = &\alpha(x) \E_{\pi} (\Delta
  S_n \mid X_n = x) \\
  & &\qquad {} + \alpha(x) \Bigl(\E_L (\Delta S_n \mid X_n = x) -
  \E_{\pi} (\Delta S_n \mid X_n = x) \Bigr) \\
  & < &\sum_k V_k \sum_{j \in C_k} x_j +\frac{1}{2} \sum_{\calN(C_0)}
  \lambda_i \hat{w}_i + \sum_{C_0} w_j \mu_j 1_{\{x_j > 0\}}  
\end{eqnarray*}
and as $V_k < 0$ for each $k$ and $x_j \geq 0$ for each $j$ the
process $S(X_n)$ is a good supermartingale under JLW at 
all but a finite subset of $\N^N$. Now Theorem \ref{Foster}(ii) implies
that the jump chain $X_n$ is positive recurrent and the ergodicity of
$X(t)$ under JLW now follows from boundedness of the event rates
$\alpha(x)$ as before. \qed  

\subsection{Proof of Theorem~\ref{local}}
In this result we are comparing behaviour of the process under
policies defined with different sets of weights so we explicitly
mention dependence upon $w$ in this section. 

Suppose that for some set of weights $w$ we have $V_1(w) < 0$. Denote
by $\hat{\pi}(w)$ a static policy that achieves drift rates $V_k(w)$ on
clusters $C_k(w)$. Using the observation in Remark \ref{weights} we
see that the processes $X(t)$ and $w'X(t)$ (for any positive weights
$w'$) are also positive recurrent under static policy
$\hat{\pi}(w)$. Now consider a policy $\hat{\pi}(w')$ that achieves
the hierarchical minimax rates for weights $w'$. By definition
\[ V_1(w') = \max_{j \in C_0} w_j^\prime \Bigl( \sum_{i \in \calN(C_0)}
\lambda_i \hat{\pi}(w' ; i)_j - \mu_j \Bigr) \leq \max_{j \in C_0}
w_j^\prime \Bigl( \sum_{i \in \calN(C_0)} \lambda_i \hat{\pi}(w; i)_j -
\mu_j \Bigr) < 0 
\]
and so $X(t)$ is positive recurrent under the static policy
$\hat{\pi}(w')$ and hence, by Theorem \ref{stability}, also under JLW
with weights $w'$.    \qed

\subsection{Preliminaries for Theorems~\ref{speeds} \& \ref{shape}}
The results of these Theorems are for the random walk $\xi(t)$ which
is obtained from the queue length process $X(t)$ by not reflecting it
at $0$.  The first result in this section is a calculation that helps
us deduce that JLW ensures that all stations in a single cluster have
the same drift rate of weighted queue length. 

The overall event rate at all states is $\alpha =
\sum_{\calN(C_0)} \lambda_i + \sum_{C_0} \mu_j$. Recall that for
the process reduced onto $D = \cup_k^K C_l$ by any static
policy $\pi \in \Pi_K$ we have to merge neighbourhoods in collections
$\sigma_i(D)$ and sum the relevant flow rates to get total flows
$\lambda_i(D)$ so, under $\pi \in \Pi_k$, the event rate at $j \in
C_k$ is $\alpha_j = \mu_j + \sum_{i \in \calN_D(D)} \lambda_i(D) \pi(i)_j$. 

We say that $x$ is \emph{properly clustered} if for each cluster $C_k$
and each $j \in C_k$, $w_j x_j < w_l x_l$ for each $l \in \cup_1^{k-1}
C_n$ and $w_j x_j > w_l x_l$ for each $l \in \cup_{k+1}^K C_n$. At a
properly clustered $x$ the event rate at station $j \in C_k$ under JLW
is $\alpha_j(x) = \mu_j + \sum_{i \in \calN_D(D) } \lambda_i(D) 1_{j \in
  \calB_i(x)}/|\calB_i(x)|$ with $D = \cup_k^K C_l$ as for $\pi \in
\Pi_k$ and $\calB_i(x) = \{ j \in S_i : w_j x_j = \min_{l \in S_i} w_l
x_l \}$ is the set of JLW routing choices for an $S_i$ arrival when
the system state is $x$.  

As before it is convenient to work with the jump chain, this time
$\xi_n$ for the random walk $\xi(t)$. For each cluster $C_k$ we will
study the process $F_k(\xi_n)$ where $F_k$ is the quadratic function  
\begin{equation}  \label{Fdef}
  F_k(x) = \frac{1}{4} \sum_{l\in C_k} \sum_{r \in C_k} \frac{(w_l
  x_l - w_rx_r)^2}{w_l w_r} =  \frac{1}{2} x^T Qx 
\end{equation}
where $Q_{lr} = -1$ for $l \neq r$ and $Q_{rr} = w_r \sum_{l \neq
  r} 1/w_l$. We write $\Delta F_k(n) = F_k(\xi_{n+1}) - F_k(\xi_n)$
and $\gamma_k = \sum_{r \in C_k} 1/w_r\,$. 

\begin{lem}  \label{dispersion}
  Consider the embedded chain $w\xi_n$ and the process $F_k(\xi_n)$. Then 
\begin{enumerate}[(i)]
\item for any $\pi \in \Pi_k$, 
\[ \alpha \E_{\pi}\Bigl( \Delta F_k(n) \mid \xi_n = x \Bigr)
  = \frac{1}{2}\sum_{j \in C_k} \alpha_j (\gamma_k w_j - 1)\,;
\]
\item for any properly clustered state $x$ 
\[  \alpha \E_L \Bigl( \Delta F_k(n) \mid \xi_n = x \Bigr) \leq
  \frac{1}{2}\sum_{j \in C_k} \alpha_j(x)(\gamma_k w_j - 1)\ .   \]
\end{enumerate}
\end{lem}

\paragraph{Proof of Lemma \ref{dispersion}} \emph{(i)} The effect of 
any policy $\pi \in \Pi_K$ is to produce 
independent random walks at stations $j \in C_k$ with the same drift
rate $w_j (\sum_{i \in \calN_D(C_k)} \lambda_i(D) \pi(i)_j - \mu_j) =
V_k$ so the calculation for this part is very 
similar to that for a zero drift random walk with independent
components. 
  
To re-use (\ref{delta_q}) we first calculate $e_j^T Qx = \gamma_k w_j x_j
- \sum_{r \in C_k} x_r$. Then observe that
$Qg = 0$ for the vector $g = (1/w_r)_{r \in C_k}$ so $F_k$ is constant
in this direction and we can translate any given $x$ in direction $g$ so that
$\sum_{C_k} x_j = 0$. After such a translation we have
$e_j^T Qx = \gamma_k w_j x_j$ when $\sum_{C_k} x_r = 0$ and also $Q_{jj}
= \gamma_k w_j - 1$. Taking
expectation under $\pi$ of (\ref{delta_q}) we have   
\begin{eqnarray*}
 \alpha \E_{\pi}\Bigl( \Delta F_k(n) \mid \xi_n = x \Bigr) 
  & = &\gamma_k \sum_{j \in C_k} w_j x_j \Bigl(\sum_{i\in
    \calN_D(C_k)} \lambda_i(D) \pi(i)_j - \mu_j\Bigr) + \frac{1}{2} \,
  \sum_{j \in C_k} \alpha_j (\gamma_k w_j - 1) \\  
  & = &\gamma_k V_k \sum_{j \in C_k} x_j + \frac{1}{2} \, \sum_{j \in
    C_k} \alpha_j (\gamma_k w_j - 1) 
\end{eqnarray*}
 and as $\sum_{j \in C_k} x_j = 0$ we have established part (i). \medskip 

\noindent\emph{(ii)} We now consider $\E_L(\Delta F_k(n) \mid \xi_n =
x)$.  As $x$ is properly clustered JLW will only route arrivals at
neighbourhoods in $\calN_D(C_k)$ in the $D$-reduced 
system to stations in $C_k$. As in the proof of Theorem
\ref{stability}, but now with $e_j^T Qx = \gamma_k w_j x_j$ when
$\sum_{C_k} x_r = 0$, 
\begin{eqnarray*}
  \lefteqn{\alpha \Bigl[ \E_L\bigl(\Delta F_k(n) \mid \xi_n = x \bigr)  -
    \E_\pi\bigl(\Delta F_k(n) \mid \xi_n = x \bigr) \Bigr] } \qquad \\
  & = &\gamma_k \sum_{j \in C_k} w_j x_j \Bigl[ \sum_{i \in \calN_D(C_k) :
    j \in S_i} 
  \lambda_i(D)  \Bigl( \frac{1_{j \in \calB_i(x)}}{|\calB_i(x)|} - \pi(i)_j
  \Bigr) \Bigr] + \frac{1}{2} \sum_{j \in C_k} (\alpha_j(x) -
  \alpha_j)(\gamma_k w_j - 1)  \\ 
  & = &\gamma_k \sum_{i \in \calN_D(C_k)} \lambda_i(D) \Bigl(
  \underline{wx}_i - \sum_{j \in S_i} \pi(i)_j w_j x_j \Bigr) +
  \frac{1}{2} \sum_{j \in C_k} (\alpha_j(x) - \alpha_j)(\gamma_k w_j - 1)  
\end{eqnarray*}
since $\sum_{j \in S_i} 1_{j \in \calB_i(x)} = |\calB_i(x)|$.  Combining
this with the result of part (i)
\begin{eqnarray*}
  \alpha \E_L\bigl(\Delta F_k(n) \mid \xi_n = x \bigr) & = &\gamma_k
  \sum_{i \in \calN_D(C_k)} \lambda_i(D) \Bigl( \underline{wx}_i - \sum_{j
    \in S_i} \pi(i)_j w_j x_j \Bigr) + \frac{1}{2} \sum_{j \in C_k}
  \alpha_j(x) (\gamma_k w_j - 1)  \\
  & \leq &\frac{1}{2} \sum_{j \in C_k} \alpha_j(x)(\gamma_k w_j - 1) 
\end{eqnarray*}
as $\underline{wx}_i \leq \sum_{j \in S_i} \pi(i)_j w_j x_j$ for each
$S_i$ at any $x$ and for any $\pi$.   \qed 

The next lemma is used in the proof of Theorem \ref{speeds} when we
show that clusters separate apart under JLW. 

\begin{lem} \label{monotonicity} 
  Consider the $D$-reduced random walk $\xi$ restricted to cluster $C$
  under JLW routing and (i) let $\zeta^+(n)$ denote the walk when
  additional arrivals (not from $\calN_D(C)$) must be routed to stations
  in $C$; (ii) let $\zeta^-(n)$ denote the walk when some arrivals to
  $\calN_D(C)$ are routed elsewhere. Suppose $\zeta^+(0) = \zeta^-(0) =
  \xi(0)$. Then $\zeta^-(n) \leq \xi(n) \leq \zeta^+(n)$ for all $n$.
\end{lem}

\paragraph{Proof of Lemma \ref{monotonicity}} 
By $x \leq z$ we mean $x_j \leq z_j$ for each $j$ and as all the $w_j
> 0$ we have $wx \leq wz$ equivalent to $x \leq z$ for all $x$, $z \in
\Z^N$. We make no assumptions about any additional arrivals or
arrivals routed elsewhere except measurability of the overall process. 

Suppose $\zeta^+(0) = \xi(0)$ and construct $\zeta^+(n)$ from the same
down jumps and arrival events as $\xi$ together with the additional
arrivals. Let $J(i, x)$ denote the station chosen by JLW for an $S_i$
arrival in state $x$. We must couple the routing processes also. In
particular if station $j = J(i, x)$ is chosen at stage $n$ for $\xi$
and $\zeta^+(n)_j = x_j$ then the same station must be chosen for
$\zeta^+$. This will work as long as $\zeta^+(n) \geq \xi(n)$ for
every $n \geq 0$ which we now show by induction.

Departures affect each process identically so cannot change
order.  At stage $n$ write $x = \xi(n) \leq z = \zeta^+(n)$. At any
additional $S_i$ arrival set $\xi(n+1) = x$ and $\zeta^+(n+1) = z
+ e_j$ where $j = J(i, z)$. At a standard $S_i$ arrival, if $J(i, x) =
J(i, z) = j$ then $\xi(n+1) = x + e_j \leq z + e_j = \zeta^+(n+1)$. If
$J(i, x) = j \neq l = J(i, z)$ then $w_j x_j \leq w_l x_l \leq w_l
z_l$ and $w_j x_j < w_j z_j$ (if $w_j z_j = w_j x_j$ the coupling
above forces $J(i, z) = j$) and again $x + e_j \leq 
z + e_l$. Hence by induction $\xi(n) \leq \zeta^+(n)$ for all $n$. 

The argument showing that $\zeta^-(n) \leq \xi(n)$ is essentially the
same but the routing coupling required is that if station $j = J(i,
z)$ is chosen at stage $n$ for $\zeta^-$ and $\xi(n)_j = z_j$ then the
same station must be chosen for $\xi$. \qed

\begin{rmk}
  This result does not extend to the queue process $X$ because the
  departure process for $X$ is dependent upon the arrival process due
  to the emptying of queues. \qquad${}_\Box$
\end{rmk}

\subsection{Proof of Theorem~\ref{speeds}}
We establish the result for a single cluster using Lemma
\ref{dispersion} and an 
inequality which we state next.  Then we use Lemma \ref{monotonicity}
to extend it to successively larger numbers of clusters.

The following generalization of Kolmogorov's maximal inequality is
Lemma 3.1 in \cite{MVW}. 
\begin{lem}\label{max_ineq}
   Let $(Y_t)_{t \in \Z^+}$ be a stochastic process on $[0, \infty)$
   adapted to a filtration $(\mathcal{F}_t)_{t \in \Z^+}$ (for example
   a function of a Markov chain). Suppose that $Y_0 = y_0$ and for
   some $b \in (0, \infty)$ and all $t \in \Z^+$ 
\[ \E\bigl( Y_{t+1} - Y_t \mid \mathcal{F}_t \bigr) \leq b \quad
a.s. \]
   Then for any $x > 0$ and any positive $t \in \Z^+$ 
\[ \P\left( \max_{0 \leq s \leq t} Y_s \geq x \right) \leq \frac{y_0 +
  bt}{x}\ .   \]
\end{lem}

Now we continue with the proof of Theorem \ref{speeds}. Suppose the
hierarchical minimax static policy results in a single cluster. We
show first that for the embedded chain $|w_j \xi_{nj} - w_l\xi_{nl}| <
n^{1-\eps}$ eventually along any sample path and we deduce the result
from this.

We again use the Lyapunov function $F_1$ introduced in
(\ref{Fdef}). We have $\max_{j, l \in C_1} | w_jx_j - w_lx_l |^2 \leq 
2\hat{w}^2F_1(x)$ at any state $x$, where $\hat{w} = \max_{C_1} w_j$.
Additionally $\max_{n/2 < r \leq n} F_1(\xi_r)/r^2 \leq 4 \max_{n/2 <
  r \leq n} F_1(\xi_r)/n^2$ along any sample path as $F_1(x) \geq 0$.
Let
\[ A_m = \Bigl\{ \max_{2^{m-1} < r \leq 2^m} \max_{j, l \in C_1} \left|
  \frac{w_j\xi_{rj} - w_l\xi_{rl}}{r} \right| \geq 2^{-m\eps}
  \Bigr\}\ . 
\]
Thus 
\[ \P_L(A_m) \leq \P_L \Bigl(\, \max_{2^{m-1} < r \leq 2^m}
  \frac{2\hat{w}^2 F_1(\xi_r)}{2^{2m}} \geq \frac{2^{-2m\eps}}{4} \,\Bigr) =
  \P_L \Bigl( \max_{2^{m-1} < r \leq 2^m} F_1(\xi_r) \geq
  \frac{2^{2m(1-\eps)}}{8 \hat{w}^2}
  \Bigr)\ .  \] 
We know from Lemma \ref{dispersion} that under JLW $\E_L(\Delta F_1(n)
\mid \xi(n) = x) \leq \gamma$ where $\gamma > 0$ is a constant. By
applying Lemma \ref{max_ineq} starting from $\xi_0 = x_0$ we now have   
\begin{eqnarray*}
   P_L (A_m) & \leq & \P_L \Bigl( \max_{0 < r \leq 2^m} F_1(\xi_r) \geq
  \frac{2^{2m(1-\eps)}}{8 \hat{w}^2}  \Bigr)\\
   & \leq &\frac{8 \hat{w}^2 \bigl( F_1(x_0) + \gamma 2^m \bigr)}{2^{2m(1-\eps)}} 
 \end{eqnarray*}
 and for $\epsilon < 1/2$ Borel-Cantelli implies only finitely many of
the $A_m$ occur. For any $n \geq 2$ we have  
\[ \left\{ \max_{r \geq n} \max_{j, l \in C_k} \left| \frac{
   w_j\xi_{rj} - w_l\xi_{rl})}{r} \right| \geq r^{-\eps}
  \right\} \subset \bigcap_{m \geq \log_2{n}} A_m 
\]
which means that for $\eps < 1/2$ there exists $n_0(\eps)$ (random)
such that $|w_j\xi_{nj} - w_l\xi_{nl}| < n^{1-2\eps}$ for all $n \geq
n_0$.  

Returning to the continuous time process $\xi(t)$ its event rate is
bounded, see (\ref{event_rate}), so this result for the jump chain
implies that all components $w_j \xi_j(t)$, $j \in C_1$, eventually have
the same drift rate.  \bigskip

Now we show that if all weighted queues in a cluster, $C_k$ say, have
the same drift rate it must be the rate obtained under the
hierarchical minimax policy. Choose $T$ large enough that 
$\max_{j, l   \in C_k} |w_l \xi_l(0) - w_j \xi_j(0)|$ is small
compared to $T$ and consider any policy $\pi$ that achieves 
\[ \max_{j, l   \in C_k} |w_l \xi_l(T) - w_j \xi_j(T)| < T^{1-\eps} \]
for some small $\eps > 0$. For each $j \in C_k$ this implies there
exist constants $V$, $\beta_j$ with $|\beta_j| < 1$ for each $j$ such
that $w_j (\xi_j(T) - \xi_j(0)) = VT + \beta_j T^{1-\eps}$. Dividing
through by $w_j T$ and summing over $j \in C_k$ we have
\[ \frac{1}{T} \, \sum_{j \in C_k} \bigl( \xi_j(T) - \xi_j(0) \bigr) =
V \sum_{j \in C_k} \frac{1}{w_j} + T^{-\eps} \sum_{j \in C_k}
\frac{\beta_j}{w_j} \ .    \] 
For large $T$ the left hand side is approximately $\sum_{i \in
  \calN_D(C_k)} \lambda_i(D) - \sum_{j \in C_k} \mu_j$ while the right
hand side is approximately $V \sum_{j \in C_k} 1/w_j$ and hence $V \to
V_k$ as $T \to \infty$ by Remark \ref{V_k}.  \bigskip

It remains to show that JLW eventually separates the clusters from any
starting configuration. We start by considering systems where
hierarchical minimax routing produces two clusters $C_1$ and
$C_2$. The optimal static policies route all arrivals at
neighbourhoods $S_i \notin \calN(C_1)$ to stations in $C_2$ while
arrivals at $S_i \in \calN(C_1)$ must be routed (by any policy) to
stations in $C_1$. The only cluster level routing error JLW can make
is to route some arrivals at $S_i \notin \calN(C_1)$ into $C_1$.

Now we employ Lemma \ref{monotonicity}.  This tells us that all
weighted queues in cluster $C_1$ eventually have speed at least
$V_1$. Also, while the cluster structure on $C_2$ may be totally
changed by the lost arrivals, no weighted queue there has speed
greater than $V_2 < V_1$ and so there exists a finite (random) time $t_0$
such that for all $t > t_0$, $w_j\xi_j(t) > w_l\xi_l(t)$ for every pair $j
\in C_1$, $l \in C_2$. For $t > t_0$ the process $w\xi(t)$ occupies
properly clustered states and so JLW no longer makes cluster level
routing errors. Now the results for single clusters imply that for
each $j \in C_k$, $w_j \xi_j(t)$ has asymptotic drift rate $V_k$ for
$k = 1$, $2$. 

Now suppose that we have established the result for systems with $K$
(hierarchical minimax) clusters and consider a system with $K + 1$
clusters. As above we see that routing errors by JLW relating to
cluster $C_1$ only act to send additional arrivals to $C_1$ and so the
weighted queue at each $j \in C_1$ eventually has speed at least
$V_1$.  The system that remains after removing $C_1$ has $K$ clusters
and initially may lose some arrivals so the maximal drift of any
weighted queue is bounded above by $V_2$.  As in the two cluster case
JLW separates $C_1$ from the rest of the system after a finite time
and the result follows by induction. \qed

\subsection{Proof of Theorem~\ref{shape}}
If $|C| = 1$ there is nothing to do so we suppose $|C| \geq 2$. We
consider a bonded sub-cluster $C \subseteq C_k$. This means that 
there is a $\pi \in \Pi_k$ such that for any stations $j$, $m \in C$
there is a path from $j$ to $m$ in the graph $G(\pi)$ with nodes
$\calN_D(C) \cup C$ and edges $\{ (S_i, j) : S_i \in \calN_D(C), j \in
C, \pi(i)_j > 0 \}$. As $C$ is finite there exists $\eps > 0$ such
that $\pi(i)_j \geq \eps$ along any such path. Similarly there exists
$\lambda^-$ such that $\lambda_i(D) \geq
\lambda^-$ for each $S_i \in \calN_D(C)$.

We modify the quadratic used in Lemma \ref{dispersion} by restricting
it to $C$ i.e.\ we use  
\[ F_C(x) = \frac{1}{4} \sum_{l, r \in C} \frac{(w_l x_l - w_r
  x_r)^2}{w_j w_r}  \ .  \] 
Repeating the calculations from Lemma
\ref{dispersion}(ii) we have  
\begin{eqnarray*}
  \alpha \E_L\bigl(\Delta F_C(n) \mid \xi_n = x \bigr) & = &\gamma_C
  \sum_{i \in \calN_D(C)} \lambda_i(D) \Bigl( \underline{wx}_i - \sum_{j
    \in S_i} \pi(i)_j w_j x_j \Bigr)    \\
  & &\qquad\qquad {} + \frac{1}{2} \sum_{j \in C} \alpha_j(x)
  (\gamma_C w_j - 1) 
\end{eqnarray*}
where $\gamma_C = \sum_{j \in C} 1/w_j$. The event rates $\alpha_j(x)$
are bounded uniformly in $x$ so 
\[ \frac{1}{2} \sum_{j \in C} \alpha_j(x)(\gamma_C w_j - 1) \leq A \]
for some constant $A$. Let $\hat{w} = \max_{j \in C} w_j$ and note
that if $F_C(x) > M^2 |C|^2 \hat{w}^2$ then $w_l x_l - w_r x_r > M$
for some pair of stations $l$, $r \in C$. 

Suppose that $w_l x_l - w_r x_r > M$ for some pair of stations $l$, $r
\in C$. As $C$ is bonded there is a loop-free path from $l$ to $r$ in
the bipartite graph $G(\pi)$. Paths in $G(\pi)$ have their nodes
alternately in $C$ and $\calN_D(C)$ and there must exist a consecutive
triple $(j, S_i, j')$ such that $j$, $j' \in S_i$ and
$w_{j'}x_{j'} - w_j x_j > M/(|C|-1)$. Thus
\[  \alpha \E_L\bigl(\Delta F_C(n) \mid \xi_n = x \bigr) \leq
\frac{-\gamma_C \lambda^- \eps M}{|C| - 1} + A   \]
which is negative for $M > A(|C| - 1)/\gamma_C \lambda^- \eps$ and
hence the process $F_C(\xi_n)$ is positive recurrent by Theorem
\ref{Foster}(ii).  \qed

\section*{Acknowledgements}
M.V. is grateful to CNPq (grants 301455/2009--0 and 472431/2009--9)
and FAPESP (thematic grant 09/52379--8),
M.M. is grateful to CNPq (grant 450787/2008--7) for partial support.

\end{document}